\RequirePackage{fix-cm}
\documentclass[smallextended]{svjour3}       
\smartqed  
\usepackage{graphicx}
\begin{document}

\title{ Discrimination of Graph Isomorphism Classes by Continuous Spectrum and Split Technique
}


\author{Ameneh Farhadian       
}


\institute{  Department of Mathematical Sciences, Sharif University of
Technology, P. O. Box 11155-9415, Tehran, I. R. Iran  \\
              \email{ a\_farhadian@mehr.sharif.ir}           
}

\date{2016-11-5}

\maketitle

\begin{abstract}

The graph isomorphism problem is a main problem which has numerous applications in different fields. Thus, finding an efficient and easy to implement method to discriminate non-isomorphic graphs is valuable. In this paper, a new method is introduced which is very simple and easy to implement, but very efficient in discriminating non-isomorphic graphs, in practice. This method does not need any heuristic attempt and based on the eigenvalues of a new matrix representation for graphs. It, almost always, separates non-isomorphic $n$-vertex graphs in time $O(n^3)$ and in worst cases such as strongly regular graphs, in time $O(n^4)$. Here, we show that this method, successfully, characterizes the isomorphism classes of studied instances of strongly regular graphs (up to 64 vertices). Strongly regular graphs are believed to be hard cases of the graph isomorphism problem.

\end{abstract}

\section {Introduction}
Graph isomorphism problem, or exact graph matching, which has numerous applications,  is the problem of determining whether  two given graphs are isomorphic. That is, there is a bijection between vertices of two graphs which preserves adjacency and non-adjacency of vertices. In order to solve this problem, one could consider each of the $n!$ bijections between vertices of two graphs, and check whether it is an isomorphism between the two graphs. If they are not isomorphic, one would need to check all $n!$ bijections to realize this fact. But, even for relatively small values of $n$, the number $n!$ is unmanageably large. 

Hence, finding an effective method with tolerable computational complexity to distinguish the isomorphism  or non-isomorphism of two graphs is very valuable. Whilst polynomial-time isomorphism-testing algorithms have been found for special classes of graphs, such as trees\cite{kelly1957congruence}, planner graphs\cite{hopcroft1974linear}, graphs of bounded degree\cite{LUKS198242} and ordered graphs\cite{Jiang19991273}, but  no polynomial runtime algorithm is known for this problem. That is, in the worst case the computational complexity of any available algorithm is exponential in the number of graph vertices \cite{johnson2005np,babai2008isomorhism}. 

One approach to solve this problem is using graph spectrum, i.e. the eigenvalues of matrix representation of a graph, as a graph  invariant.  Since, the matrix representation of a simple graph is symmetric, graph spectrum is computable in time $O(n^3)$. Thus, a graph spectrum can be regarded as a fast graph invariant to distinguish graph isomorphism classes, provided that it has been chosen properly. From this point of view, a good matrix representation for graphs is a matrix that produces minimum possible number of cospectral graphs. This idea widely has been used and investigated in several papers, such as spectral determination \cite{cvetkovic2010spectral}, spectral characterization  \cite{MR3003429},  spectrum as a graph invariant \cite{MR2212507},  comparing graphs via graph spectrum \cite{wilson2008study}, graph clustering \cite{cvetkovic2011graph}.

Until now, several alternative matrix representations for graphs have been proposed in the literatures, such as adjacency matrix, Laplacian matrix, signless Laplacian \cite{MR2070541,MR2022290}, normalized Laplacian \cite{MR1421568} and heat kernel \cite{Bai:2005:MEG:2145086.2145089}.
The spectrum of all these representations may be used to characterize graphs. Wilson and Zue\cite{wilson2008study}  have compared these spectrums according to cospectrality. They have shown that there are some graphs which even the combination of these spectrums is not able to separate them. In other words, even the combination of these spectra is not sufficient to determine the graph isomorphism classes, uniquely. To avoid cospectral graphs, Emms et al.\cite{MR2212507} have introduced a matrix representation inspired by the notion of coined quantum walk. They have shown the graph spectrum with respect to this matrix is successful to separate many known families of cospectral graphs. But, there are some non-isomorphic graphs, like a 14-vertex graphs pair, which are not distinguishable by their spectrum. In addition, the size of the matrix that they have proposed is in the order of $O(n^2)$. 

Here, we proposes a new matrix representation for graphs. The eigenvalues of this matrix, as a graph invariant, provides a simple polynomial-time algorithm to efficiently solve the graph isomorphism problem in practice.  Continuous spectrum and Split technique which are defined in this paper, mutually, provide a fast computable graph invariant to solve some computationally hard problems, such as graph isomorphism and graph automorphism, in practice. 

The continuity matrix  representation for graphs  is defined in the next section. The eigenvalues of this matrix which is in terms of a variable is called continuous spectrum. This spectrum is equipped with split method in the third section to make it more efficient in separating non-isomorphic graphs. The paper is concluded in the fourth section. Some derived relations and theorems about continuous spectrum are given in the Appendix.

\section {Continuous Spectrum}
 In this section, we define continuity matrix representation of graphs. Also, we describe why we expect the eigenvalues of such matrix to behave well in separating non-isomorphic graphs.\\
At first, let us review some used definitions and notations. 
 The adjacency matrix of a graph $G$ is a matrix $A = [a_{ij} ]$ where $a_{ij}=1$  if $v_{i}$ is adjacent to $v_{j}$ and $a_{ij}=0$ otherwise.
The Laplacian matrix of a graph is $L=D-A$. Matrix $D$ is diagonal degree matrix, whose diagonal elements are given by the vertex degrees $D(u,u)=d_{u}$.  The sign-less  Laplacian  and the normalized Laplacian matrix, respectively, are $ |L|=D+A$  and $ \widehat{L}=D^{-1/2}LD^{1/2}$.\\
For a given spectrum, two graphs are cospectral if they share the same spectrum.
 A graph $G$ is said to be determined by its spectrum, or DS\cite{MR2022290}, if any graph $H$ that is cospectral with
$G$ is also isomorphic to $G$. 
\begin{definition}
Let $A$ be the adjacency matrix of an $n$-vertex graph $G$ and $k$ be a fixed non-negative integer.  We define continuity matrix $$C( \alpha):=B-Diag(B) $$ where $ B=  A+\alpha A^{2}+ \cdots + \alpha^{k} A^{k+1}$ and  $| \alpha | <1$ is a real variable. 
 $Diag(B)$ is a diagonal matrix with the diagonal entries of matrix $B$. Continuous spectrum of $G$ is the eigenvalues of $C(\alpha)$.

\end{definition}

Clearly, the matrix $C(\alpha)$ is the same adjacency matrix, if $k=0$ or $\alpha=0$. Continuous spectrum is the eigenvalues of matrix $C(\alpha)$ and, consequently, is in terms of $\alpha$. To compute the continuous spectrum for a given graph, it is sufficient to assign an arbitrary value to $\alpha$.  

\begin{figure}[ht]
\centerline{\includegraphics[width=4cm]{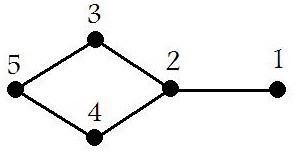}}
\caption{\label{graf}\small A sample graph and its continuity matrix }
\end{figure}

Computer results show that small value of $k$ is enough to separate non-isomorphic graphs by the continuous spectrum. For instance, taking $k=2$ is enough to discriminate graph isomorphism classes for graphs with up to $9$ vertices. 
As an example, the continuity matrix representation of the graph Fig.\ref{graf} is given( for $k=2$), here.

$$C(\alpha)= \left[ \begin{array}{c c c c c}
0 & 1+3 \alpha^2 & \alpha & \alpha & 2 \alpha^2 \\
1+3 \alpha^2 & 0 & 1+5\alpha^2 & 1+5\alpha^2 & 2 \alpha \\
\alpha & 1+5 \alpha^2 & 0 & 2\alpha & 1+4\alpha^2 \\
\alpha & 1+5 \alpha^2 & 2\alpha & 0 & 1+4\alpha^2\\
2\alpha^2 & 2 \alpha &  1+4\alpha^2&  1+4\alpha^2 & 0
\end{array} \right]$$
\begin{figure}[ht]
\centerline{\includegraphics[width=12cm]{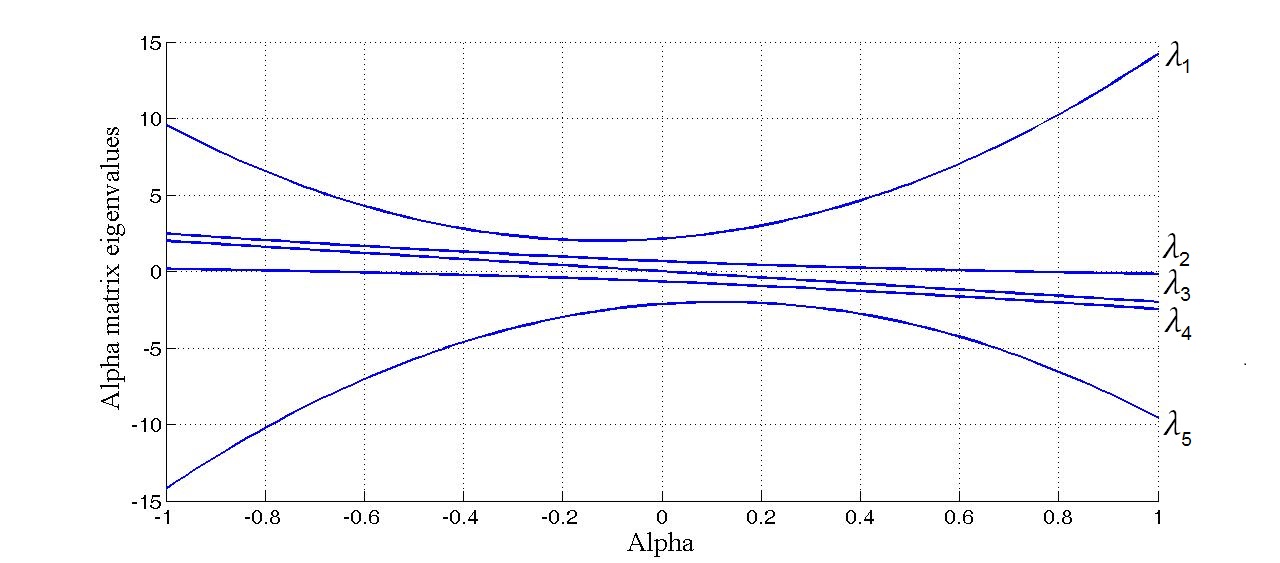}}
\caption{\label{diagram}\small Continuous spectrum digram in terms of $\alpha$ }
\end{figure}
 Figure \ref{diagram} depicts the continuous spectrum of the mentioned graph in terms of $\alpha$. Since this graph is bipartite, its diagram is symmetric about the origin. The proof of this property for bipartite graphs is given in the Appendix.

The aim of defining such spectrum is to enhance the ability of graph spectrum in distinguishing non-isomorphic graphs. We would like  a graph spectrum to play the role of a fast computable invariant to discriminate non- isomorphic graphs.

Let us see why we expect continuous spectrum to be efficient in discriminating non-isomorphic graphs. Let $A$ be the adjacency matrix of a graph $G$ and $\lambda$ is an eigenvalue of $A$ with the corresponding eigenvector $ u $.
Thus, we have $ Au=\lambda u $. In other words,
$$ \left\lbrace \begin{array}{c}
a_{1,1} u_{1}+a_{1,2} u_{2}+ \cdots +a_{1,n} u_{n}=\lambda u_{1} \\
a_{2,1} u_{1}+a_{2,2} u_{2}+ \cdots +a_{2,n} u_{n}=\lambda u_{2} \\
\vdots \\
a_{n,1} u_{1}+a_{n,2} u_{2}+ \cdots +a_{n,n} u_{n}=\lambda u_{n}
\end{array} \right. $$
The above system of equations says that the summation of the values of neighbors of a vertex should be the value of that vertex multiplied by $ \lambda$.
In this view, $\lambda $ is a possible value under this condition that summation of neighbors values of each vertex should be the value of that vertex multiplied by $\lambda$. Thus, only the neighbor vertices are regarded and other vertices are ignored in each equation.\\
A natural question arises here. Why only the direct neighbors are considered and other vertices which are two or more steps away are ignored? How can we consider the effect of vertices which are not neighbors of vertex $v_{i}$ in the $i$-th equation of the above system of equations while preserving the importance of the neighbors?\\

Therefore, we want to modify adjacency matrix to take into account the effect of non-neighbor vertices while preserving the importance of the neighbors. We use an attenuator factor, namely $\alpha$, to reduce the effect of far vertices. Hence, in each equation we enter the effect of vertices which are $ i+1$ steps away, by adding their corresponding entries multiplied by $\alpha ^i$. Also, if there are more than one way between $v_{i}$ and $v_{j}$, we enumerate all of them. Exactly speaking, we consider $k$-length walks and their numbers. We know that  $(i,j)$-entry of $A^{k}$ indicates the  number of $k$-length walks between $ v_{i} $ and $v_{j} $. Thus, we take into account  the effect of vertex $v_{j}$ which is $k$ edges away from $v_i$ by  $(i,j)$-entry of  $A^{k}$  multiplied by $\alpha^{k-1}$. Then, we put zero on the diagonal. Therefore, the desired matrix is what you see in the above definition.

Graphs which are determined by their spectrum, abbreviated by DS\cite{MR2022290}, are those that do not have any cospectral mate. In order to compare the successfulness of one spectrum in characterizing a set of graphs, we count the number of graphs which are determined by that spectrum.
The distinguishing ability of adjacency spectrum, Laplacian , signless Laplacian and normalized Laplacian spectrum are compared with continuous spectrum  in Table \ref{small}.
\begin{table}[h]\label{smallgraphs}
\begin{center}
\begin{tabular}{c|c |c| c|c|c|c  }
    n &  \# & A & L & $\vert L \vert$ & $L_{N} $& C$(\alpha)^{\alpha =.7}_{  k=2}$\\
\hline
 6 &156& 146& 152 &140 &142 &  156\\
 7  &1044&934& 914& 942 & 992 &  1044 \\
 8 & 12346&10624& 10579 &11145& 12145& 12346 \\
 9 & 274668&223827&  270732 & 256262&274590& 274668\\
\end{tabular}
\caption { \label{small}\small{The number of DS graphs for small graphs with respect to different matrix representations}}
\end{center}
\end{table}
According to Table \ref{smallgraphs}, as we expected due to  the definition of continuity matrix, its eigenvalues is  more efficient in discriminating non-isomorphic small graphs compared with other matrices. Since continuity matrix is symmetric, the complexity of computing continuous spectrum for a definite value of $\alpha$  is $O(n^3)$ for a graph with $n$ vertices.

\section{Split Spectrum}
In the previous section, continuous spectrum for separating graph isomorphism classes of graphs was defined. In spite of its efficiency in separating non-isomorphic graphs, there are  non-isomorphic graphs that the continuous spectrum does not discriminate them. Here, Split technique for indexing graph vertices is defined. It provides a proficient graph  isomorphism class distinguisher.

Let $A_1$ and $A_2$ be the adjacency matrices of two non-isomorphic graphs $G_1$ and $G_2$. If for any integer $i$, $Diag(A_1^i)=Diag(A_2^i)=kI$ , then the continuous spectrum is identical for both of them and, consequently, is not able to separate them. For instance, the diagonal of $A^t$ for strongly regular graphs with the same parameters, is the same and equal to $k.I$. Thus, continuous spectrum cannot distinguish them. Here, a new labeling for vertices is proposed which uses a graph spectrum. Combination of the following technique with continuous spectrum  provides an efficient invariant, even for characterizing strongly regular graphs. 

\begin{definition}
The split spectrum on a vertex $v$  of a graph is $$ S_v:= (S_{N_v} ,S_{ \backslash{N_v}}) $$ where $S_{N_v}$ and $S_{\backslash {N_v}}$ are, respectively, the spectrum of the induced subgraph on the neighbors and non-neighbors of the vertex $v$.\\ The split spectrum of a graph $G$ is
$$ Split  Spec(G):=\lbrace Sv \vert v \in V(G) \rbrace $$

 \end{definition}
The split spectrum can be computed with respect to any graph spectrum, such as adjacency spectrum, Laplacian spectrum or continuous spectrum. The computation results show that split spectrum with respect to continuous spectrum, i.e. Split continuous spectrum, is a powerful distinguisher for graph isomorphism classes. To examine the efficiency of split continuous spectrum, we have investigated it for some set of graphs.
 Since continuous spectrum is the same adjacency spectrum when $\alpha=0$, the continuous spectrum inherits the successfulness of adjacency spectrum in  discrimination of graphs which do not have any cospectral graphs. Thus, it is sufficient to check its efficiency for discrimination of cospectral graphs (with respect to adjacency matrix). \\
Strongly regular graphs provide relatively large sets of cospectral graphs with respect to $A$ and also to $L$, $\vert L \vert$,  $L_{N}$ and even $C(\alpha)$. Therefore, they are good and challenging choices to check the efficiency of split continuous spectrum.  Furthermore, the family of strongly regular graphs has long been identified as a hard case for the graph isomorphism problem \cite{read1977graph} and the best existing GI algorithms for them are exponential \cite{babai2013faster,babai1980complexity,spielman1996faster}. Thus, finding a simple and fast polynomial time algorithm to practically solve the graph isomorphism problem for strongly regular graphs is valuable. \\
A $k$-regular graph with $n$ vertices is strongly regular SRG($n$,$k$,$r$,$s$), if every two adjacent vertices have $r$ common neighbors and every two non-adjacent vertices have $s$ common neighbors \cite{bose1963strongly}. In Table \ref{srg}, some sets of  strongly regular graphs and their numbers are given. They are obtained from \cite{mckay,ted}.\\
\begin{table}
\begin{tabular}{c c  c c c c c c }
    SRG(n,k,r,s) &\# &   Split Sp. $A$ & Split Sp. $C(\alpha)_{k=3}^{\alpha =.3}$ & Split Sp. $C(\alpha)_{k=10}^{\alpha =.3,.7}$  \\
\hline
 (26,10,3,4)  & 10& .20&1&1\\
 (29,14,6,7)  &41&.02&1&1\\
 (35,16,6,8) &3854&.90&1&1\\
 (35,18,9,9) & 227 &.55&1&1\\
 (36,14,4,6) &180&.20&.92&1 \\ 
 (37,18,8,9) &6760& .00&1&1\\
 (40,12,2,4) &28& .04&.79&1\\
 (45,12,3,3) &78& .77 &1&1\\ 
(50,21,8,9) &18&1 &1&1\\ 
(64,18,2,6) &167&.01 &.72&1\\
\end{tabular}
\caption{\label{srg}\small{The fractions of DS graphs for SRGs with respect to split spectrum   }}
\end{table}
We have measured the distinguishing ability of split spectrum with respect to different graph spectra  by the fraction of DS graphs that they produce in Table \ref{srg}. We know that the eigenvalues of all known graph representing matrices, such as adjacency matrix,  Laplacian  matrix and even continuity matrix is the same for strongly regular graphs with the same parameters. But, as computer results show split continuous spectrum  is able to easily distinguishes the isomorphism classes of such graphs.\\
 In Table \ref{srg}, the results of using split spectrum for discriminating graph isomorphism classes of SRGs are shown. This  table shows that the split spectrum with respect to adjacency matrix is not able to successfully separate non-isomorphic classes of SRG graphs. In contrast, the combination of continuous spectrum with split technique provides an efficient tool to  distinguish graph isomorphism classes. As we see in the table,  taking $k=10$ and $\alpha= .3$ , $.7$ is enough to discriminate all mentioned instances of strongly regular graphs with $26$ to $64$ vertices. 

Since the computational complexity of this method is $O(n^4)$  for $n$-vertex graphs, it significantly outperforms to other methods which are exponential. Thus, it is a valuable approach to solve graph isomorphism problem in practice, even for strongly regular graphs.
\section{Conclusion}
The eigenvalue of a matrix representation of graphs, i.e. graph spectrum, is invariant to arrangement of graph vertices and is computable in $O(n^3)$.  Therefore, it is a good invariant to separate non-isomorphic graphs. But, usual spectrums such as adjacency spectrum, Laplacian spectrum and even combination of them suffer from existing cospectral graphs. For example, for all strongly regular graphs with the same parameter, both adjacency and Laplacian matrix produce the same eigenvalues.

In this paper, a new matrix representation for graphs, i.e. continuity matrix, is proposed. This representation has improved the ability of graph spectrum in recognizing non-isomorphic graphs. The definition of continuity matrix is based on the extension of the equations that eigenvalues satisfy in. In a graph adjacency matrix, each equation is established for a vertex and its neighbors, while in the new defined matrix each equation is written for all vertices. This fact makes the graph spectrum more engaged in graph structure and, consequently, enhances its ability in discrimination of non-isomorphic graphs. Therefore, finding cospectral graphs with respect to this new matrix is difficult.

For discrimination of  some special graphs such as strongly regular graphs, which have been known as hard cases of the isomorphism problem, continuous spectrum is not sufficient, alone. Thus, split spectrum was proposed to improve continuous spectrum in separating non-isomorphic graphs. We have shown by computations that split continuous spectrum, which  is computable in $O(n^4)$, successfully characterizes  the isomorphism classes of many known families of strongly regular graphs up to 64 vertices.

Continuous spectrum and split method, mutually,  provide an efficient invariant to discriminate non-isomorphic graphs.
Such graph spectrum which is highly entangled with graph structure benefits us in solving some computationally hard problems, such as graph isomorphism  and graph automorphism problem, in practice. Split continuous spectrum, which does not need any heuristic attempt,  provides a simple, efficient and fast graph invariant to distinguish graph isomorphism classes in  average time $ O(n^3)$ and $O(n^4)$ in worst case.

\section*{Acknowledgment}
I would like to thank Ted Spence and Richard Wilson for providing graph data.
\bibliographystyle{plain}

\section{Appendix: Some results about Continuous Spectrum}
The graph shown in Fig. \ref{graf} is a bipartite graph and the digram of its continuous spectrum, as we see in Fig. \ref{diagram}  is symmetric about the origin. In other words, in bipartite graphs the eigenvalues of $C(\alpha) \vert _{\alpha=\alpha_{0}}$  and  $C(\alpha) \vert _{\alpha=-\alpha_{0}}$  for any arbitrary $\alpha_{0} $ is the same, but with a negative sign. The next theorem proves this fact.

\begin{theorem}
Let $G$ be a bipartite graph. Then
$$ \left( \begin{array}{c c c c}
\lambda_{1}(\alpha) & \lambda_{2}(\alpha) & \cdots & \lambda_{s}(\alpha) \\
m_{1} & m_{2} & \cdots & m_{s}
\end{array} \right) =\left( \begin{array}{c c c c}
-\lambda_{s}(-\alpha) & -\lambda_{s-1}(-\alpha) & \cdots & -\lambda_{1}(-\alpha) \\
m_{s} & m_{s-1} & \cdots & m_{1}
\end{array} \right)$$

The above notation is from  \cite{biggs:1993} and displays the eigenvalues and their multiplicity.
\end{theorem}
\textbf{Proof}: It is sufficient to show that if $ \lambda $ is an eigenvalue of $ C(\alpha) $ then $ -\lambda $ is an eigenvalue of $ C(-\alpha) $ and vice versa.
Assume the continuity matrix of $ G $ is $ C(\alpha)= \left( \begin{array}{c c }
X(\alpha) & Z(\alpha)  \\
 Z^{T}(\alpha)  & Y(\alpha)\\ \end{array} \right) $ where  $ X(\alpha) $ and $ Y(\alpha) $ are square matrices corresponding to each part of bipartite graph. Clearly, these two matrices are odd functions of $ \alpha $. And the functions $ Z(\alpha)$ and $Z^{T}(\alpha)$ are even functions of $ \alpha $.\\
Let the partitioned vector $ u=(u_{X},u_{Y}) $ be an eigenvector
of $ C(\alpha)$ with eigenvalue $ \lambda $ where $ u_{X}$ and $u_{Y} $respectively relate to parts  $ X $ and $ Y $. Now, $ u'=(u_{X},-u_{Y}) $ is an eigenvector of $ A(-\alpha)$  with eigenvalue  $ -\lambda $, because  \\
$ C(-\alpha)u'= \left( \begin{array}{c c}
X(-\alpha) & Z(-\alpha)  \\
 Z^{T}(-\alpha)  & Y(-\alpha)\\ \end{array} \right)u'=\left( \begin{array}{c c}
-X(\alpha) & Z(\alpha)  \\
 Z^{T}(\alpha)  & -Y(\alpha)\\ \end{array} \right)\left( \begin{array}{c  c}
u_{X}   \\
-u_{Y}\\ \end{array} \right)=\left( \begin{array}{c  c}
-\lambda u_{X}   \\
\lambda u_{Y}\\ \end{array} \right)= -\lambda u'$. \\
It follows that $ \lambda $ and $ -\lambda $ are respectively eigenvalues of  $ C(\alpha)$ and $ A(-\alpha)$ with the same multiplicity. $\diamond$

 The corollary of above theorem is a familiar theorem in algebraic graph theory.
\begin{corollary}
If graph $G$ is bipartite, its adjacency spectrum is symmetric about the origin \cite{MR1829620}.
\end{corollary}
\textbf{Proof}: It is sufficient to take $ \alpha=0 $ in the previous theorem. \\

There are several theorems for ordinary graph spectrum which are remained also true for continuous spectrum. Some of them are given here. Their proof is similar to what is given in algebraic graph theory books for adjacency spectrum, such as \cite{biggs:1993}.
\begin{theorem}
Let $C(\alpha)$ be the continuity matrix of graph $G$, and $\pi$ a permutation of $V(G)$. Then $\pi$ is an automorphism of $G$ if and only if $PC(\alpha)=C(\alpha)P$ , where $P$ is the permutation matrix representing $\pi $.
\end{theorem}
\begin{theorem}
Let $ \lambda (\alpha) $ be a simple eigenvalue of $G$ for a specified $ \alpha $, and let $x(\alpha )$ be a corresponding eigenvector(with real components). If the permutation matrix $P$ represents an automorphism of $G$ then $P \times x(\alpha)=\pm x(\alpha)$.
\end{theorem}
\begin{theorem}
If there is an $\alpha$ that eigenvalues of $C(\alpha) \vert_{\alpha}$ are simple, every automorphism of $G$ (apart form the identity) has order 2.
\end{theorem}

\end{document}